\documentstyle[amssymb]{aipproc}
\newtheorem{theorem}{\sc Theorem}
\newtheorem{definition}{\sc Definition}
\righthead{\sc Deformation Quantization}
\lefthead{\sc Daniel Sternheimer}
\begin{document}
\title{Deformation Quantization: \\ Twenty Years After\footnote
{{\it This review is dedicated to the memory of our good friend Ryszard
R\c aczka}, in whose honor the meeting in \L \'od\'z ``Particles, Fields
and Gravitation", was held in April 1998. The e-print is infinitesimally
updated from the August 1998 version to be published in the Proceedings
entitled {\it Particles, Fields and Gravitation}, edited by Jakub
Rembieli\'nski for AIP Press.}}
\author{Daniel Sternheimer}
\address{Laboratoire Gevrey de Math\'ematique Physique, CNRS \ ESA 5029\\
D\'epartement de Math\'ematiques, Universit\'e de Bourgogne\\
BP 400, F-21011 Dijon Cedex France. \\ {\rm e-mail:}
 {\tt dastern@u-bourgogne.fr}}
\date{}
\maketitle
\vspace{-10mm}

\begin{abstract}
We first review the historical developments, both in physics and in
mathematics, that preceded (and in some sense provided the background of)
deformation quantization. Then we describe the birth of the
latter theory and its evolution in the past twenty years, insisting on the
main conceptual developments and keeping here as much as possible on the
physical side. For the physical part the accent is put on its relations
to, and relevance for, ``conventional" physics. For the mathematical part
we concentrate on the questions of existence and equivalence,
including most recent developments for general Poisson manifolds; we touch
also noncommutative geometry and index theorems, and relations with
group theory, including quantum groups. An extensive (though very incomplete)
bibliography is appended and includes background mathematical literature.
\end{abstract}
\vspace{-6mm}

\section*{I. Background}
In this Section we briefly present the fertile ground which was needed in
order for deformation quantization to develop, even if from an abstract
point of view one could have imagined it on the basis of Hamiltonian
classical mechanics. Indeed there are two sides to ``deformation
quantization". The philosophy underlying the r\^ole of {\it deformations
in physics} has been consistently put forward by Flato since more than
30 years and was eventually expressed by him in \cite{Fl82} (see also
\cite{Fa90,Fl98}). In short, the passage from one level of physical theory to
another, more refined, can be understood (and might even have been predicted)
using what mathematicians call deformation theory. For instance one
passes from Newtonian physics to special relativity by deforming the
invariance group (the Galilei group
${\rm{SO}}(3)\cdot{\Bbb R}^3\cdot{\Bbb R}^4$) to the Poincar\'e group
${\rm{SO}}(3,1)\cdot{\Bbb R}^4$ with deformation parameter
$c^{-1}$, where $c$ is the velocity of light. There are many other examples
among which {\it quantization} is perhaps the most seminal.
\newpage

As a matter of fact it seems that the idea that quantum mechanics is some
kind of deformed classical mechanics has been, almost from the beginning
of quantum theory, ``in the back of the mind of many physicists" (after we
came out with the preprint of \cite{BFFLS78}, a scientist even demanded
that we quote him for that!). This is attested by the notion of classical limit
and even more by that of semi-classical approximation, a good presentation
of which can be found in \cite{Vo77}. But the idea remained hidden ``in the
back of the minds" for a long time, in particular due to the apparently
insurmountable ``quantum jump" in the nature of observables -- and probably
also because the mathematical notion of deformation and the relevant
cohomologies were not available. A long maturation was needed which
 eventually gave birth to full-fledged deformation quantization about
20 years ago \cite{BFFLS78}.

A word of caution may be needed here. It is possible to intellectually imagine
new physical theories by deforming existing ones. Even if the mathematical
concept associated with an existing theory is mathematically rigid, it may
be possible to find a wider context in which nontrivial deformations exist.
For instance the Poincar\'e group may be deformed to the simple (and therefore
rigid in the category of Lie groups) anti De Sitter group ${\rm{SO}}(3,2)$ very
popular recently  -- though it had been studied extensively by us 15-20 years
ago \cite{AFFS81}, resulting in particular in a formulation of {\it QED with
photons dynamically composed of two singletons} \cite{FF88} in AdS universe.
As is now well-known, there exist deformations of the Hopf algebras
associated with a simple Lie group (these ``quantum groups" \cite{Dr86},
which are in fact an example of deformation quantization \cite{BFGP94},
have been extensively studied and applied to physics). Nevertheless such
intellectual constructs, even if they are beautiful mathematical theories,
need to be somehow confronted with physical reality in order to be taken
seriously in physics. So some physical intuition is still needed when using
deformation theory in physics.

\subsection*{I.1 Weyl Quantization and Related Developments}
We assume the reader somewhat familiar with (classical and) quantum
mechanics. In an ``impressionist" fashion we only mention the names of
Planck, Einstein, de~Broglie, Heisenberg, Schr\"odinger and finally
Hermann Weyl. What we call here ``deformation quantization" is related to
Weyl's quantization procedure. In the latter \cite{We31}, starting with
a classical observable $u(p,q)$, some function on phase space
${\Bbb R}^{2\ell}$ (with $p, q \in {\Bbb R}^{\ell}$), one associates an
operator (the corresponding quantum observable) $\Omega(u)$ in the
Hilbert space $L^2({\Bbb R}^{\ell})$ by the following general recipe:
\begin{equation} \label{weyl}
 u \mapsto \Omega_w(u) = \int_{\Bbb R}^{2\ell} \tilde{u}(\xi,\eta)
{\exp}(i(P.\xi + Q.\eta)/\hbar)w(\xi,\eta)\ d^\ell \xi d^\ell \eta
\end{equation}
where $\tilde{u}$ is the inverse Fourier transform of $u$, $P_\alpha$ and
$Q_\alpha$ are operators satisfying the canonical commutation relations
$[P_\alpha , Q_\beta] = i\hbar\delta_{\alpha\beta}$
($\alpha, \beta = 1,...,\ell$), $w$ is a weight function and the integral
is taken in the weak operator topology.  What is now called normal ordering
corresponds to choosing the weight
$w(\xi,\eta) = {\exp}(-{1\over 4}(\xi^2 \pm \eta^2))$,
standard ordering (the case of the usual pseudodifferential operators in
mathematics) to $w(\xi,\eta) = {\exp}(-{i\over 2}\xi\eta)$ and the original
Weyl (symmetric) ordering to $w = 1$.
An inverse formula was found shortly afterwards by Eugene Wigner \cite{Wi32}
and maps an operator into what mathematicians call its symbol by a kind
of trace formula. For example $\Omega_1$ defines an isomorphism of Hilbert
spaces between $L^2({\Bbb R}^{2\ell})$ and Hilbert-Schmidt operators on
$L^2({\Bbb R}^{\ell})$ with inverse given by
\begin{equation} \label{EPW}
u=(2\pi\hbar)^{-\ell}\, {\rm{Tr}}[\Omega_1(u)\exp((\xi.P+\eta.Q)/i\hbar)]
\end{equation}
and if $\Omega_1(u)$ is of trace class one has
${\rm{Tr}}(\Omega_1(u))=(2\pi\hbar)^{-\ell}\int u \, \omega^\ell$ where
$\omega^\ell$ is the (symplectic) volume $dx$ on ${\Bbb R}^{2\ell}$.
Numerous developments followed in the direction of phase-space methods, many
of which can be found described in \cite{AW70}. Of particular interest to us
here is the question of finding an interpretation to the classical function
$u$, symbol of the quantum operator $\Omega_1(u)$; this was the problem posed
(around 15 years after \cite{Wi32}) by Blackett to his student
Moyal. The (somewhat na\"{\i}ve) idea to interpret it as a probability density
had of course to be rejected (because $u$ has no reason to be positive)
but, looking for a direct expression for the symbol of a quantum commutator,
Moyal found \cite{Mo49} what is now called the Moyal bracket:
\begin{equation}
M(u,v) = \nu^{-1} \sinh(\nu P)(u,v) =
P(u,v) + \sum^\infty_{r=1}\nu^{2r}P^{2r+1} (u,v)   \label{Moyal}
\end{equation}
where $2\nu=i\hbar$, $P^r(u,v)=\Lambda^{i_1j_1}\ldots
\Lambda^{i_rj_r}(\partial_{i_1\ldots i_r}u)(\partial_{j_1\ldots j_r} v)$
is the $r^{\rm{th}}$ power ($r\geq1$) of the Poisson bracket
bidifferential operator $P$, $i_k, j_k = 1,\ldots,2\ell$, $k=1,\ldots,r$
and $(\Lambda^{i_kj_k}) = {0\,-I\choose I\,0}$.
To fix ideas we may assume here $u,v\in C^\infty({\Bbb R}^{2\ell})$ and the
sum taken as a formal series (the definition and convergence for various
families of functions $u$ and $v$ was also studied, including in
\cite{BFFLS78}). A similar formula for the symbol of a product
$\Omega_1(u)\Omega_1(v)$ had been found a little earlier \cite{Gr46} and
can now be written more clearly as a (Moyal) {\it star product}:
\begin{equation} \label{star}
u \ast_M v = \exp(\nu P)(u,v) = uv + \sum^\infty_{r=1}\nu^{r}P^{r}(u,v).
\end{equation}
Several integral formulas for the star product have been introduced and
the Wigner image of various families of operators (including bounded
operators on $L^2({\Bbb R}^{\ell})$) were studied, mostly after deformation
quantization was developed (see e.g. \cite{Da80,Ma86,Ka86}).
An adaptation to Weyl ordering of the mathematical notion of
pseudodifferential operators (ordered, like differential operators, 
``first $q$, then $p$") was done in~\cite{GLS68} -- and the converse
in \cite{Ho79}.
Starting from field theory, where normal (Wick) ordering is essential
(the r\^ole of $q$ and $p$ above is played by $q\pm ip$), Berezin
\cite{Be74,BS70} developed in the mid-seventies an extensive study of what
he called ``quantization", based on the correspondence principle and Wick
symbols. It is essentially based on K\"ahler manifolds and related to
pseudodifferential operators in the complex domain \cite{BG81}.
However in his theory (which we noticed rather late), as in the
studies of various orderings \cite{AW70}, the important concepts of
{\it deformation} and {\it autonomous} formulation of quantum mechanics
in general phase space are absent.

\subsection*{I.2 Classical Mechanics on General Phase Space \\
 and its Quantization}
Initially classical mechanics, in Lagrangean or Hamiltonian form, assumed
implicitly a ``flat" phase space ${\Bbb R}^{2\ell}$, or at least considered
only an open connected set thereof. Eventually more general configurations
were needed and so the mathematical notion of manifold, on which mechanics
imposed some structure, was needed. This has lead in particular to using the
notions of symplectic and later of Poisson manifolds, which have been
introduced also for purely mathematical reasons. One of these reasons has
to do with families of infinite-dimensional Lie algebras, which date back
to works by \'Elie Cartan at the beginning of this century and regained a
lot of popularity (including in physics) in the past 30 years.

A typical example can be found with Dirac constraints \cite{Di64}: second
class Dirac constraints restrict phase space from some ${\Bbb R}^{2\ell}$
to a symplectic manifold $W$ imbedded in it (with induced symplectic form),
while first class constraints further restrict to a Poisson manifold
with symplectic foliation (see e.g. \cite{FLS76}). Some of the references
where one can find detailed information on the symplectic approach to
classical (Hamilton) mechanics are \cite{LM87,AM78,GS84} and (which
includes the derivation of symplectic manifolds from Lagrangean mechanics)
\cite{So97}. The question of quantization on such manifolds was certainly
treated by many authors (including in \cite{Di64}) but did not go beyond
giving some (often useful) recipes and hoping for the best.

A first systematic attempt started around 1970 with what was called soon
afterwards {\it geometric quantization} \cite{Ko70}, a by-product of Lie group
representations theory where it gave significant results \cite{AK71,Ki76}.
It turns out that it is geometric all right, but its scope as far as
quantization is concerned has been rather limited since few classical
observables could be quantized, except in situations which amount
essentially to the Weyl case considered above.
In a nutshell one considers phase-spaces $W$ which are coadjoint orbits of
some Lie groups (the Weyl case corresponds to the Heisenberg group with the
canonical commutation relations ${\frak{h}}_\ell$ as Lie algebra); there one
defines a ``prequantization" on the Hilbert space $L^2(W)$ and tries
to halve the number of degrees of freedom by using polarizations (often
complex ones, which is not an innocent operation as far as physics is
concerned) to get a Lagrangean submanifold ${\cal L}$ of dimension half
that of $W$ and quantized observables as operators in $L^2({\cal L})$;
``Moyal quantization" on a symplectic groupoid 
${\Bbb{R}}_1^{2\ell}\times {\Bbb{R}}_2^{2\ell}$ was obtained therefrom
in \cite{GV95}. A recent exposition can be found in \cite{Wo92}.

\subsection*{I.3 Pseudodifferential Operators and Index Theorems}
One may argue that physicists had invented the theory of distributions
(with Dirac's $\delta$) and symbols of pseudodifferential operators (with
standard ordering) much before mathematicians developed the corresponding
theories. These may also be considered as belonging to the large family
of examples of a fruitful interaction between physics and mathematics, even
if in the latter case (symbols) it seems that the two developments were
largely independent at the beginning, and in fact converged only with the
advent of deformation quantization.

In this connection one should not forget that there is a significant
difference in attitude to Science (with notable exceptions):
in physics a very good idea may be enough to earn you a Nobel prize
(with a little bit of luck, enough PR and provided you live long enough
to see it well recognized); in mathematics one usually needs to have,
young enough, several good ideas and prove that they are really good
(often with hard work, because ``problems worthy of attack prove their worth
 by hitting back") in order to be seriously considered for the Fields medal.
This rather ancient difference may explain Goethe's sentence
(much before Nobel and Fields):
``Mathematicians are like Frenchmen: They translate everything into their
own language and henceforth it is something completely different".

In the fifties \cite{CZ57} the notion of Fourier integral operators
was introduced, generalizing and making precise the sometimes heuristic
calculus of ``differential operators of noninteger order". Soon it
evolved into what are now called pseudo\-differential operators, defined
on general manifolds \cite{Ho66} and as indispensable to theories
of partial differential equations as distributions (with which they are
strongly mixed). But what gained to this tool fame and respectability
among all mathematicians was the proof, in 1963 (and following years for
various generalizations) of the so-called {\it index theorem} for elliptic
(pseudo)differential operators on manifolds by Atiyah, Singer, Bott, Patodi
and others \cite{AS63,Pa65,CS63,Gi74}. The (analytical) index of a linear
map between two vector spaces is defined (when both terms are finite) as the
dimension of its kernel minus the codimension of its image. Elliptic partial
differential operators $d$ on compact manifolds $X$ have an index $i(d)$,
equal (and this is the original theorem) to a ``topological index"
which depends only on topological invariants associated with the manifold
(the Todd class $\tau(X)$ of the complexified cotangent bundle of $X$ and
the fundamental class $[X]$) and on a cohomological invariant ${\rm{ch}}d$
(a Chern character) associated with the symbol of the principal part of the
operator $d$. To give the flavor of the result we write a precise formula
(see e.g. Atiyah's lecture in \cite{CS63}), valid for compact manifolds
(without or with boundary):
\begin{equation} \label{ASI}
i(d)=\langle({\rm{ch}}d)\tau(X),[X]\rangle.
\end{equation}
The existence of such a formula had also been conjectured by Gel'fand.
The proof is very elaborate and one cannot avoid doing it also for
pseudodifferential operators. Topological arguments and factorization
(which imposes consideration of continuous symbols -- this was my share
in \cite{CS63}) permit eventually to reduce the proof of the equality to
the cases of the Dirac operator on even dimensional manifolds and one
particular (convolution) operator on the circle.
There have been numerous developments in a wide array of mathematical
domains provoked by this seminal result. The formula itself has been
very much generalized, including to ``algebraic" index theorems where the
algebra of pseudodifferential operators is replaced by more abstract algebras
(this is an major ingredient in noncommutative geometry \cite{Co94} and
is strongly related to star products \cite{CFS92,NT95}).
Throughout the theory a capital r\^ole is played by the symbol $\sigma(d)$
(the classical function associated with the standard-ordered operator $d$).
Note that the principal part of a differential operator is independent of
the ordering, but eventually the whole symbol was used. In the proof of the
reduction one needs an expression for the symbol of a product of operators,
given by an integral formula analogous to a star product. So mathematicians
had been using star products (albeit corresponding to a different ordering and
without formal series development in some parameter like $\hbar$)
before they were systematically defined. This permitted eventually to give
original proofs of existence of star products on quite general manifolds
\cite{LB95,Gu95} by adapting techniques and results developed \cite{BG81}
in the theory of pseudodifferential operators.

\subsection*{I.4 Cohomologies and Deformation Theory}
In an often ignored section of a paper, I.E. Segal \cite{Se51} and
(independently) a little later Wigner and Inon\"u \cite{IW53,Sa61}
have introduced in the early fifties a kind of inverse \cite{LN67}
to the mathematical notion of {\it deformation} of Lie groups and algebras,
notion which was precisely defined only in 1964 by
Murray Gerstenhaber \cite{Ge64}. That inverse was called {\it contraction}
and typical examples (mentioned at the beginning of this Section)
are the passage from De Sitter to Poincar\'e groups (by taking the limit
of zero curvature in space-time) or from Poincar\'e to Galilei (by taking
$c^{-1}\rightarrow 0$). Intuitively speaking a contraction
is performed by neglecting in symmetries, at some level of physical reality,
a constant (like $c^{-1}$) which has negligible impact at this level but
significant effects at a more ``refined" level. Note that this may be
realized mathematically in varying generality (e.g. \cite{Sa61} is more general
than \cite{IW53} but both have for inverse a Gerstenhaber deformation).
The notion of deformation of algebras, which may be seen as an
outcome of the notion (introduced a few years before) of deformations of
complex analytic structures \cite{KS58}, gives rise to a better defined
mathematical theory which, for completeness, we shall briefly present
in the following two subsections. All this is by now well-known, even to
physicists, and we shall keep details to a minimum, referring the reader
interested in more details to papers and textbooks cited here and
references quoted therein.

It also turns out that recently we had to introduce deformations which are
even more general than those introduced by Gerstenhaber
(see \cite{DFST97,Pi97,Na98} and \cite{Fl98} in these Proceedings) in order
e.g. to quantize Nambu mechanics. They are still inverse to some
contraction procedure, applied (like deformation quantization, where the
algebra is that of classical observables and the parameter Planck's constant)
to algebras which are not geometrical symmetries.
So one should keep an open mathematical mind, let physics be a guide and
develop if needed completely new mathematical tools.
This is true {\it physical mathematics}, in contradistinction with standard
mathematical physics where one mainly applies existing tools or with
theoretical physics where mathematical rigor is too often left aside and a
good physical intuition (which Dirac certainly had e.g. when he
worked with his $\delta$ ``function") is then required as a guide.

\subsubsection*{I.4.1 Hochschild and Chevalley-Eilenberg cohomologies}
Let first $A$ be an {\it associative} algebra (over some commutative ring
${\Bbb K}$) and for simplicity we consider it as a module over itself with
the adjoint action (algebra multiplication); the generalization to cohomology
valued in a general module is straightforward.
A {\it $p$-cochain} is a $p$-linear map $C$ from $A^p$ into (the module) $A$
and its {\it coboundary} $bC$ is given by
\begin{eqnarray}
bC(u_0,\ldots,u_p)&=&u_0C(u_1,\ldots,u_p)-C(u_0u_1,u_2,\ldots,u_p)+
\cdots\nonumber\\
&+&(-1)^pC(u_0,u_1,\ldots,u_{p-1}u_p)+(-1)^{p+1}C(u_0,\ldots,u_{p-1})u_p.
\end{eqnarray}
One checks that we have here what is called a complex, i.e. $b^2=0$.
We say that a $p$-cochain $C$ is a {\it $p$-cocycle} if $bC=0$. We denote
by ${\cal Z}^p(A,A)$ the space of $p$-cocycles and by ${\cal B}^p(A,A)$ the
space of those $p$-cocycles which are coboundaries (of a $(p-1)$-cochain).
The $p$th {\it Hochschild cohomology} space (of $A$ valued in $A$) is
defined as ${\cal H}^p(A,A)={\cal Z}^p(A,A)/{\cal B}^p(A,A)$. {\it Cyclic
cohomology} is defined using a bicomplex which includes the Hochschild
complex and we shall briefly present it at the end of this review in the
example of interest for us here.

For {\it Lie algebras} (with bracket $\{\cdot,\cdot\}$) one has a similar
definition, due to Chevalley and Eilenberg \cite{CE48}. The $p$-cochains
are here skew-symmetric, i.e. linear maps $B: \wedge^p A\longrightarrow A$,
and the Chevalley coboundary operator $\partial$ is defined on a $p$-cochain
$B$ by (where ${\hat{u_j}}$ means that $u_j$ has to be omitted):
\begin{eqnarray}
\partial C(u_0,\ldots,u_p)=
&\sum_{j=0}^p&(-1)^j\{u_j,C(u_0,\ldots,{\hat{u_j}},\ldots,u_p)\}\nonumber\\
&+&\sum_{i<j}(-1)^{i+j}C(\{u_i,u_j\},u_0,\ldots,{\hat{u_i}},\ldots,
{\hat{u_j}},\ldots,u_p).
\end{eqnarray}
Again one has a complex ($\partial^2=0$), cocycles and coboundaries spaces
$Z^p$ and $B^p$ (resp.) and by quotient the {\it Chevalley cohomology} spaces
$H^p(A,A)$, or in short $H^p(A)$; the collection of all cohomology spaces
is often denoted $H^*$.

\subsubsection*{I.4.2 Gerstenhaber theory of deformations of algebras}
Let $A$ be an algebra. By this we mean an {\it associative, Lie} or {\it Hopf}
 algebra, or a {\it bialgebra}. Whenever needed we assume it is also a
{\it topological} algebra, i.e. endowed with a locally convex topology for
which all needed algebraic laws are continuous. For simplicity we may think
that the base (commutative) ring ${\Bbb K}$ is the field of complex
numbers ${\Bbb C}$ or that of the real numbers ${\Bbb R}$. Extending it to
the ring ${\Bbb K}[[\nu]]$ of formal series in some parameter $\nu$
gives the module ${\tilde A} = A[[\nu]]$, on which we can consider the
preceding various algebraic (and topological) structures.
\medskip

{\it I.4.2.1 Deformations and cohomologies}.
A concise formulation of a Gerstenhaber deformation of an algebra (which
we shall call in short a {\it DrG-deformation} whenever a confusion may
arise with more general deformations) is \cite{Ge64,GS88,BFGP94}:
\begin{definition}
A deformation of such an algebra $A$ is a ${\Bbb K}[[\nu]]$-algebra
${\tilde A}$ such that ${\tilde A}/\nu {\tilde A} \approx A$.
Two deformations ${\tilde A}$ and ${\tilde A'}$ are said equivalent if they
are isomorphic over ${\Bbb K}[[\nu]]$ and ${\tilde A}$ is said trivial if
it is isomorphic to the original algebra $A$ considered by base field
extension as a ${\Bbb K}[[\nu]]$-algebra.
\end{definition}
Whenever we consider a topology on $A$, ${\tilde A}$ is supposed to be
topologically free. For associative (resp. Lie algebra) Definition 1 tells us
that there exists a new product $\ast$ (resp. bracket $[\cdot,\cdot]$) such
that the new (deformed) algebra is again associative (resp. Lie).
Denoting the original composition laws by ordinary product 
(resp. $\{\cdot,\cdot\}$) this means that, for $u,v\in A$
(we can extend this to $A[[\nu]]$ by ${\Bbb K}[[\nu]]$-linearity) we have:
\begin{eqnarray}
u\ast v &=& uv + \sum_{r=1}^\infty \nu^r C_r(u,v) \label{a}\\
\left[u, v\right] &=& \{u,v\} + \sum_{r=1}^\infty \nu^r B_r(u,v)\label{l}
\end{eqnarray}
where the $C_r$ are Hochschild 2-cochains and the $B_r$ (skew-symmetric)
Chevalley 2-cochains, such that for $u,v,w\in A$ we have
$(u\ast v) \ast w=u\ast (v\ast w)$ and ${\cal{S}}[[u,v],w]=0$, where
${\cal{S}}$ denotes summation over cyclic permutations. At each level $r$
we therefore need to fulfill the equations ($j,k\geq1$):
\begin{eqnarray}
D_r(u,v,w)&\equiv&\sum_{j+k=r}(C_j(C_k(u,v),w)-C_j(u,C_k(v,w))=bC_r(u,v,w)
\label{adefcond}\\
E_r(u,v,w)&\equiv&\sum_{j+k=r}{\cal{S}}B_j(B_k(u,v),w)=\partial B_r(u,v,w)
\label{ldefcond}
\end{eqnarray}
where $b$ and $\partial$ denote (respectively) the Hochschild and Chevalley
coboundary operator. In particular we see that for $r=1$ the driver
$C_1$ (resp. $B_1$) must be a 2-cocycle. Furthermore, assuming one has shown
that (\ref{adefcond}) or (\ref{ldefcond}) are satisfied up to some order $r=t$,
a simple calculation shows that the left-hand sides for $r=t+1$ are then
3-cocycles, depending only on the cochains $C_k$ (resp. $B_k$) of order
$k\leq t$. If we want to extend the deformation up to order $r=t+1$
(i.e. to find the required 2-cochains $C_{t+1}$ or $B_{t+1}$), this cocycle
has to be a coboundary (the coboundary of the required cochain):
{\it The obstructions to extend a deformation from one step to the next lie
in the 3-cohomology}. In particular (and this was Vey's trick) if one can
manage to pass always through the null class in the 3-cohomology, a cocycle
can be the driver of a full-fledged (formal) deformation.

For a (topological) {\it bialgebra} (an associative algebra $A$ where we have
in addition a coproduct $\Delta : A \longrightarrow A \otimes A$ and the
obvious compatibility relations), denoting by $\otimes_\nu$ the tensor
product of ${\Bbb K}[[\nu]]$-modules, we can identify
${\tilde A}\, {\hat{\otimes}}_\nu {\tilde A}$ with
$(A\, {\hat{\otimes}}A)[[\nu]]$, where ${\hat{\otimes}}$ denotes the algebraic
tensor product completed with respect to some operator topology (e.g.
projective for Fr\'echet nuclear topology), we similarly have a deformed
coproduct ${\tilde \Delta } = \Delta + \sum_{r=1}^\infty \nu^r D_r$,
$D_r \in {\cal L}(A, A {\hat{\otimes}}A)$ and in this context appropriate
cohomologies can be introduced. Here we shall not elaborate on these, nor on
the additional requirements for Hopf algebras, referring for more details
to original papers and books; there is a huge literature on the subject,
among which we may mention \cite{Dr86,GS90,Bo92,BFGP94,BFP92,SS93,Ta90}
and references quoted therein.
\medskip

{\it I.4.2.2 Equivalence} means that there is an isomorphism
$T_\nu=I+\sum_{r=1}^\infty \nu^r T_r$, $T_r\in{\cal L}(A,A)$ so that
$T_\nu(u*'v)=(T_\nu u*T_\nu v)$ in the associative case, denoting by
$*$ (resp. $*'$) the deformed laws in ${\tilde A}$ (resp. ${\tilde A'}$);
and similarly in the Lie case.
In particular we see (for $r=1$) that a deformation is trivial at order 1 if
it starts with a 2-cocycle which is a 2-coboundary. More generally, exactly
as above, we can show \cite{BFFLS78} that if two deformations are equivalent
up to some order $t$, the condition to extend the equivalence one step
further is that a 2-cocycle (defined using the $T_k$, $k\leq t$) is the
coboundary of the required $T_{t+1}$ and therefore {\it the obstructions to
equivalence lie in the 2-cohomology}. In particular, if that space is null,
all deformations are trivial.
\medskip

{\it I.4.2.3 Unit}. An important property is that a {\it deformation of an
associative algebra with unit} (what is called a unital algebra) is
again unital, and {\it equivalent to a deformation with the same unit}.
This follows from a more general result of Gerstenhaber (for deformations
leaving unchanged a subalgebra) and a proof can be found in \cite{GS88}.
\medskip

{\it I.4.2.4}. In the case of (topological) {\it bialgebras} or {\it Hopf}
algebras, {\it equivalence} of deformations has to be understood as
an isomorphism of (topological) ${\Bbb K}[[\nu]]$-algebras, the isomorphism
starting with the identity for the degree 0 in $\nu$. A deformation is again
said {\it trivial} if it is equivalent to that obtained by base field
extension. For Hopf algebras the deformed algebras may be taken
(by equivalence) to have the same unit and counit, but in general not the
same antipode.

\subsubsection*{I.4.3 Examples of special interest: the differentiable cases}
Consider the algebra $N=C^\infty(X)$ of functions on a differentiable
manifold $X$. When we look at it as an associative algebra acting on itself
by pointwise multiplication, we can define the corresponding Hochschild
cohomologies. Now let $\Lambda$ be a skew-symmetric contravariant two-tensor
(possibly degenerate) defined on $X$, satisfying $[\Lambda,\Lambda]=0$
in the sense of the supersymmetric Schouten-Nijenhuis bracket
\cite{Ni55,Sc54} (a definition of which, both intrinsic and in terms of
local coordinates, can be found in \cite{BFFLS78,FLS74}).
Then the inner product $P(u,v)=i(\Lambda)(du\wedge dv)$ of $\Lambda$ with
the 2-form $du\wedge dv$, $u,v\in N$,  defines a
{\it Poisson bracket} $P$: it is obviously skew-symmetric,
satisfies the Jacobi identity because $[\Lambda,\Lambda]=0$ and the Leibniz
rule $P(uv,w) = P(u,w)v + uP(v,w)$. It is a bidifferential 2-cocycle for the
(general or differentiable) Hochschild cohomology of $N$, skewsymmetric
of order $(1,1)$, therefore \cite{BFFLS78} nontrivial and thus defines an
infinitesimal deformation of the pointwise product on $N$. We say that $X$,
equipped with such a $P$, is a {\it Poisson manifold} \cite{BFFLS78,Li77}.

When $\Lambda$ is everywhere nondegenerate ($X$ is then necessarily of
even dimension $2\ell$), its inverse $\omega$ is a closed everywhere
nondegenerate 2-form ($d\omega=0$ is then equivalent to $[\Lambda,\Lambda]=0$)
and we say that $X$ is {\it symplectic}; $\omega^\ell$ is a volume element
on $X$. Then one can in a consistent manner work with differentiable cocycles
\cite{BFFLS78,FLS74} and the differentiable Hochschild $p$-cohomology space
${\cal H}^p(N)$ is \cite{HKR62,Ve75} that of all skew-symmetric contravariant
$p$-tensor fields, and therefore is infinite-dimensional. Thus, except when
$X$ is of dimension 2 (because then necessarily ${\cal H}^3(N)=0$),
the obstructions belong to an infinite-dimensional space where they may be
difficult to trace. On the other hand, when $2\ell=2$, any 2-cocycle can
be the driver of a deformation of the associative algebra $N$:
``anything goes" in this case; some examples for ${\Bbb R}^2$ can be found
in \cite{Ve75}.

Now endow $N$ with a Poisson bracket: we get a Lie algebra and can look at
its Chevalley cohomology spaces. Note that $P$ is bidifferential of order
$(1,1)$ so it is important to check whether the Gerstenhaber theory is
{\it consistent} when restricted to {\it differentiable} cochains (both of
arbitrary order and of order at most 1), especially since the general
(Gelfand-Fuks) cohomology is very complicated \cite{GKS72} (but in fact the
pathology arises only when non-continuous cochains are allowed).
This is a nontrivial question and we gave it a positive answer
\cite{FLS74,BFFLS78}; in brief, if a coboundary is differentiable, it is
the coboundary of a differentiable cochain.

Again, since $P$ is of order (1,1), we first studied the
1-differentiable cohomologies. When the cochains are restricted to be of
order $(1,1)$ with no constant term (then they annihilate constant functions,
which we write ``n.c." for ``null on constants") it was found \cite{Li74}
that the Chevalley-Eilenberg cohomology $H^*_{\rm{1-diff,n.c.}}(N)$
of the Lie algebra $N$ (acting on itself with the adjoint representation)
is exactly the de Rham cohomology $H^*(X)$. Thus
${\rm{dim}}H^p_{\rm{1-diff,n.c.}}(N)=b_p(X)$, the $p$th Betti number of the
manifold $X$. Without the n.c. condition one gets a slightly more complicated
formula \cite{Li74}; in particular if $X$ is symplectic with an exact
2-form $\omega=d\alpha$, one has here
$H^p_{\rm{1-diff}}(N)=H^p(X)\oplus H^{p-1}(X)$.

All this allowed us ``three musketeers" to study in 1974 what we called in
 \cite{FLS74} {\it 1-differentiable deformations} of the Poisson bracket Lie
algebra $N$ and to give some applications \cite{FLS76}. In particular
we noticed that the ``pure" order $(1,1)$ deformations correspond to a
deformation of the 2-tensor $\Lambda$; allowing constant terms and taking
the deformed bracket in Hamilton equations instead of the original Poisson
bracket gave (at this classical level) a kind of friction term.

Shortly afterwards, triggered by our works, a ``fourth musketeer" J. Vey
\cite{Ve75} noticed that in fact
${\rm{dim}}H^2_{\rm{diff,n.c.}}(N)=1+{\rm{dim}}H^2_{\rm{1-diff,n.c.}}(N)$
and that $H^3_{\rm{diff,n.c.}}(N)$ is also finite-dimensional, which allowed
him to study differentiable deformations. Incidentally, in the
${\Bbb R}^{2\ell}$ case, Vey rediscovered (independently, because he
ignored it) the Moyal bracket. The latter was then rather ``exotic" and
few authors (except for a number of physicists, like \cite{AW70} or 
J. Pleba\'nski who described it in Polish lecture notes \cite{Pl69} he gave
us later) paid any attention to it. In {\it Mathematical Reviews} this
bracket, for which \cite{Mo49} is nowadays often quoted, is not even
mentioned in the review! We then came back to the problem \cite{FLScr76},
this time with differentiable deformations, and deformation quantization
was conceived.
\newpage

\section*{II. Deformation Quantization \\ and its Ramifications}

\subsection*{II.1 The Birth of Deformation Quantization}

Though we had mentioned the main features in 1976 and 1977 in short papers
\cite{FLS76,BFFLS77}, meetings \cite{FS80} and a long preprint, it is only
with the publication of the latter \cite{BFFLS78}, our first major
(and often quoted) contribution in this new domain, that what eventually
became known as {\it deformation quantization} \cite{Ws94} took off.
Incidentally, as for the two other parts of our ``deformation trilogy"
(see e.g. \cite{FF78,FF88} and \cite{FPS77,FS77,FS80}) which deal (resp.)
with singleton physics and with nonlinear evolution equations, true
recognition was slow to come (it took about twenty years). Let me stress
once more (and this will become evident with the section on physical
applications), that the important and most original conceptual aspect is
that {\it quantization} is here an {\it autonomous theory} based on a
{\it deformation} of the composition law of classical observables, not
on a radical change in the nature of the observables. In addition to this
important {\it conceptual advantage}, our approach is more general (simple
examples can be given); it can be shown to coincide with the conventional
(operatorial) approach in known applications (see below) whenever a (possibly
generalized) Weyl mapping can be defined; it also paves the way to better
conventional quantizations in field theory (e.g. on the infinite-dimensional
symplectic manifold of initial conditions for nonlinear evolution
equations) via a kind of cohomological renormalization.

\subsubsection*{II.1.1 Differentiable deformations and star products}
Let $X$ be a differentiable manifold (of finite, or possibly infinite,
dimension; to be precise, in the former case, we assume it is locally of
finite dimension, paracompact and Hausdorff, and by differentiable we mean
infinitely differentiable; the base field may be ${\Bbb R}$ or ${\Bbb C}$).
We assume given on $X$ a {\it Poisson structure} (a Poisson bracket $P$).
\begin{definition}
A {\it star product} is a deformation of the associative algebra of functions
$N=C^\infty(X)$ of the form $\ast = \sum_{n=0}^\infty \nu^n C_n$
where for $u,v\in N$, $C_0(u,v)=uv$, $C_1(u,v)-C_1(v,u)=2P(u,v)$ and the
$C_n$ are bidifferential operators (locally of finite order).
We say a star product is {\it strongly closed} if
$\int_X (u*v-v*u) dx =0$ where $dx$ is a volume element on $X$.
\end{definition}
\noindent {\it Remark 1. a}. The parameter $\nu$ of the deformation is in
physical applications taken to be $\nu= {\frac{i}{2}}\hbar$.

{\it b}. Using equivalence one may take $C_1=P$. The latter is the case of
Moyal, but other orderings like standard or normal do not verify this
condition (only the skew-symmetric part of $C_1$ is $P$).
Again by equivalence, in view of Gerstenhaber's result mentioned in I.4.2.3,
we may take cochains $C_r$ which are without constant term (what we called
n.c. or null on constants). In fact, in the original paper \cite{BFFLS78},
we considered only this case and we also concentrated on ``Vey products"
\cite{Li82} for which the cochains $C_r$ have the same parity as $r$ and
have $P^r$ for principal symbol in any Darboux chart, with $X$ symplectic;
when $X$ is symplectic of dimension $2\ell$ with symplectic
form $\omega$, the (Liouville) volume element is $dx=\omega^\ell$.
  
{\it c}. It is also possible to consider star products for which the
cochains $C_n$ are allowed to be slightly more general. Allowing them
to be {\it local} ($C_n(u,v)=0$ on any open set where $u$ or $v$ vanish)
gives nothing new \cite{CGW80}. (Note that this is not the same as requiring
the whole associative product to be local; in fact \cite{Ru84} the latter
condition is very restrictive and, like true pseudodifferential operators,
a star product is a nonlocal operation). In some cases (e.g. for star
representations of Lie groups) it may be practical to consider
pseudodifferential cochains. As far as the cohomologies are
concerned, it has been recently shown \cite{Gu97,PG97} that as long as one
requires at least continuity for the cochains, the theory is the same as
in the differentiable case. Incidentally this indicates that one has to
go beyond continuous cochains to get the pathological features of the
general Gelfand-Fuks cohomology of infinite Lie algebras.

Also, due to formulas like (\ref{EPW}) and the relation with Lie algebras
(see II.4.1), it is sometimes convenient \cite{NT95,Fe96,BW96} to take
${\Bbb K}[\nu^{-1},\nu]]$ (Laurent series in $\nu$, polynomial in
$\nu^{-1}$ and formal series in $\nu$) for the ring on which the deformation
is defined. Again, this will not change the theory.

{\it d}. By taking the corresponding commutator
$[u,v]_\nu=(2\nu)^{-1}(u\ast v - v \ast u)$, since the skew-symmetric
part of $C_1$ is $P$, we get a deformation of the Poisson bracket
Lie algebra $(N,P)$. This is a crucial point because
(at least in the symplectic case) we know the needed Chevalley cohomologies
and (in contradistinction with the Hochschild cohomologies) they are small
\cite{Ve75,DLG84}. The interplay between both structures gives existence
and classification; in addition it will explain why (in the symplectic case)
the classification of star products is based on the 1-differentiable
cohomologies, hence ultimately on the de Rham cohomology of the manifold.
 
\subsubsection*{II.1.2 Invariance and covariance}
Since the beginning \cite{BFFLS78} we realized that there is a big difference
between Poisson brackets and their deformations, from the point of view of
geometric invariance. Indeed, while a Poisson bracket $P$ is (by definition)
invariant under all symplectomorphisms, i.e. transformations of the manifold
$X$ which preserve the symplectic form $\omega$ (generated by the flows
$x_u=i(\Lambda)(du)$ defined by Hamiltonians $u\in N$), already on
${\Bbb R}^{2\ell}$ one sees easily that its powers $P^r$, hence also the
Moyal bracket (\ref{Moyal}), are invariant only under flows generated by
Hamiltonians $u$ which are polynomials of maximal order 2, forming the
``affine" symplectic Lie algebra
${\frak{sp}}({\Bbb R}^{2\ell})\cdot{\frak{h}_\ell}$. For other orderings
the invariance is even smaller (only ${\frak{h}_\ell}$ remains). For general
Vey products the first terms of a star product are \cite{BFFLS78,Li82}
$C_2 = P^2_\Gamma + bH$ and $C_3 = S^3_\Gamma + T +3\partial H$.
Here $H$ is a differential operator of maximal order 2, $T$ a 2-tensor
corresponding to a closed 2-form, $\partial$ the Chevalley coboundary
operator. $P^2_\Gamma$ is given (in canonical coordinates) by an expression
similar to $P^2$ in which usual derivatives are replaced by covariant
derivatives with respect to a given symplectic connection $\Gamma$
(a torsionless connection with totally skew-symmetric components when all
indices are lowered using $\Lambda$). $S^3_\Gamma$ is a very special cochain
given by an expression similar to $P^3$ in which the derivatives are
replaced by the relevant components of the Lie derivative of $\Gamma$ in
the direction of the vector field associated to the function ($u$ or $v$).
Fedosov's algorithmic construction \cite{Fe94} shows that the symplectic
connection $\Gamma$ plays a r\^ole at all orders. Therefore the invariance
group of a star product is a subgroup of the {\it finite-dimensional} group
of symplectomorphisms preserving a connection. Its Lie algebra is
${\frak{g}}_{0}=\{a\in N; [a,u]_\nu=P(a,u)\, \forall u\in N \}$,
{\it preferred observables} Hamiltonians for which the classical and
quantum evolutions coincide.
We are thus lead to look for a weaker notion and shall call a star product
{\it covariant} under a Lie algebra ${\frak{g}}$ of functions if
$[a,b]_\nu = P (a,b)$ $\forall a,b \in {\frak{g}}$. It can be shown
\cite{ACMP83} that $*$ is ${\frak{g}}$-covariant iff there exists a
representation $\tau$ of the Lie group $G$ whose Lie algebra is
${\frak{g}}$ into ${\rm{Aut}}(N[[\nu]]; *)$ such that
$\tau_g u = (Id_N + \sum^\infty_{r=1}\nu^r\tau^r_g)(g.u)$ where
$g \in G, u \in N$, $G$ acts on $N$ by the natural action induced by the
vector fields associated with ${\frak{g}}$, $(g\cdot u)(x)=u(g^{-1} x)$,
and where the $\tau^r_g$ are differential operators
on $W$. Invariance of course means that the geometric action preserves the
star product: $g\cdot u \ast g \cdot v = g\cdot (u \ast v)$. This is the
basis for the theory of star representations which we shall briefly
present below (II.4.1).

\subsection*{II.2 Existence and Classification\\ of Deformation Quantizations}

\subsubsection*{II.2.1 Symplectic finite-dimensional manifolds; reduction}
As early as 1975, Vey \cite{Ve75} had shown that on a symplectic manifold $X$
with $b_3(X)=0$, there exists a globally defined deformation of the Poisson
bracket $P$. He did this by a careful study allowing him to show (by induction)
that at each order of deformation (each order of $\hbar$), he can manage
to pass via the zero class of the finite-dimensional obstructions space
$H^3_{\rm{diff,n.c.}}(N)$. This was later easily extended to star products
\cite{NV81} and in essence tells us that we can in this case ``glue" Moyal
products defined on local charts (equivalence will take care of intersection
of two charts and the vanishing of $b_3(X)$ permits to do it in a way
compatible with multiple intersections).

The restriction $b_3(X)=0$ seemed purely technical and indeed already in
\cite{BFFLS78}, with the important case of the hydrogen atom ($X=T^*(S^3)$),
we showed that it is not essential. The latter case was generalized by Gutt
to $X=T^*(M)$ where $M$ is a Lie group \cite{Gu83} (see also regular 
star representations in \cite{CG82}) or more generally a parallelizable
manifold. Shortly afterwards De Wilde and Lecomte were able to find a proof
of existence, first for a general cotangent bundle, then for exact
symplectic manifolds and finally for a general symplectic manifold.
The latter required at first a very abstract proof which the authors
eventually made more ``palatable" \cite{DL83} along the lines used
by the Japanese group (gluing Moyal on Darboux charts),
and the question of invariant star products was also studied \cite{DLM85}.

What is behind the scene for cotangent bundles $T^*(M)$, is that there one has
globally defined ``momentum coordinates" which permits to work with globally
defined differentiations on $M$ and polynomials in them (differential
operators). From there a natural step forward is to localize everything
on a general symplectic manifold $X$ and to work with a bundle ${\cal{W}}(X)$
of {\it Weyl algebras}; a Weyl algebra $W_\ell$ is the enveloping algebra of
the Heisenberg canonical commutation relations Lie algebra ${\frak{h}}_\ell$,
possibly completed to formal series, and the product there is the
${\rm{Sp}}(2\ell)$-invariant Moyal product. The ``miracle" is that this bundle
has a flat connection and a global section; therefore locally defined (Moyal)
star products on polynomials (in any given Darboux chart) can be glued together
to a globally defined star product. (A Darboux chart on a symplectic manifold
is a chart with local coordinates the $q$'s and $p$'s of physicists, i.e.
where $\omega=\sum_{j=1}^{2\ell}dq_j\wedge dp_j$).
This line of conduct was taken by Fedosov since 1985 (in an obscure paper
\cite{Fe85} which was made detailed and precise later \cite{Fe94,Fe96}),
using a ``germ" approach (infinitesimal neighborhoods), symplectic
connections and an algorithmic construction of the flat connection on the
Weyl bundle, canonically constructed starting from a given symplectic
connection on $X$.
Independently a Japanese group \cite{OMY91} obtained also a general proof
of existence, gluing together Moyal products defined on Darboux charts by
``projecting" from the Weyl bundle, and the method could also be easily
adapted to give the existence \cite{OMY91} of {\it closed star products}
\cite{CFS92}. 
For more details on these (especially Fedosov's) constructions we refer to
the original papers and to the reviews in \cite{BCG98,GR98,Ws94}.
The relation between the ``Russian" approach and the ``Belgian" one was
made \cite{De95} by a famous Belgian mathematician (with a Russian wife),
translating both into his own language of ``gerbes".
Eventually the Fedosov construction was shown (cf. e.g. \cite{NT95}) to
be ``generic" in the sense that any differentiable star product is
equivalent to a Fedosov star product.
\medskip

Now an important tool in symplectic mechanics is that of {\it reduction}
\cite{MW74} caused by the action of an invariance group $G$ and subsequent
reduction of the algebra of observables. In fact, already in \cite{BFFLS78},
a reduction of this general type was used in connection with the hydrogen
atom. Fedosov has recently showed \cite{Fe98} that the classical reduction
theory can be ``quantized" in the same conditions, i.e. that reduction
commutes with $G$-invariant deformation quantization (at least for $G$
compact); note however that, as shown in a simple example \cite{Wa98},
reductions may give nonequivalent star products.
\medskip

\noindent{\it Remark 2: connection with 1-differentiable deformations}.
We have indicated that equivalence classes of star products are in one-to-one
correspondence with formal series in the deformation parameter $\nu$ with
coefficients in $H^2(X)$, i.e. series of the form
$[\Lambda] + \nu [H^2(X)][[\nu]]$, where $\Lambda$ is the (nondegenerate)
Poisson 2-tensor given on $X$. Since $H^2(X)$ classifies 1-differentiable
deformations of the Lie algebra $(N,P)$, one expects that {\it any Fedosov
deformation can be obtained by a sequence of successive 1-differentiable
deformations of the initial Poisson bracket}, starting with any given one,
at least for $X$ symplectic. This is indeed true.
Intuitively, one starts -- if needed by equivalence -- with a star product
n.c. at order~1, $(u,v) \mapsto uv + {\nu}P(u,v) + O(\hbar^2)$.
Associativity is satisfied mod$\nu^2$. Then one takes a 1-differentiable
infinitesimal deformation of $P$ of the form  $P+\nu P_1 +\cdots$,
corresponding in fact to an infinitesimal deformation
$\omega+\nu\omega_1+O(\nu^2)$. Fedosov tells us that we can find $P'_1$
and higher order terms so that the new product
$(u,v)\mapsto uv +{\nu}(P(u,v)+\nu P_1(u,v))+\nu^2 P'_1 (u,v) + \cdots $
is associative to order 2 in $\nu$. The classes of the 2-tensors associated
with $P_1$ give all possible choices to order 2. One does the same at the
next step with $(P + \nu P_1 + \nu^2 P_2)$ and $P'_2$, and so on. Indeed,
in \cite{Bo98}, where at every order in $\nu$ the effect of adding a
de Rham 2-cocycle was traced in the star product, Bonneau gave a detailed
proof of this fact. So by ``plugging in" 1-differentiable deformations of
Poisson brackets one can cover all possible equivalence classes, starting
from any given Fedosov star product.

In a more abstract form a similar conclusion is a consequence of results
developed in the nice review \cite{GR98}. The mathematically oriented
reader will find there a detailed presentation (in a \v{C}ech cohomology
context) of the equivalence question and related problems.

\subsubsection*{II.2.2 Poisson finite-dimensional manifolds}
As we explained earlier, physics (e.g. with first class Dirac constraints)
requires sometimes manifolds which have a Poisson bracket $P$, but a
degenerate one and are therefore not symplectic. These are called {\it Poisson
manifolds} \cite{BFFLS78} and they are foliated with symplectic leaves.
A typical example is the dual of a finite-dimensional Lie algebra, foliated
by coadjoint orbits (see e.g. \cite{Ki76}). As this example shows (even in
the case of the Heisenberg Lie algebra ${\frak{h}}_\ell$, where one must not
forget the trivial orbit), Poisson manifolds are in general not regular;
a {\it regular} Poisson manifold is foliated with symplectic leaves of
constant even dimension (before introducing Poisson manifolds we had
considered \cite{FLS74} ``canonical manifolds", regular Poisson manifolds
where the leaves have codimension 1). In the regular case the theory we just
explained extends in a straightforward manner. Some non conclusive attempts
were made in the general case, most notably in \cite{OMY94} following
more ``traditional" lines and \cite{Vn97} in the direction indicated by
Kontsevich's formality conjecture \cite{Kt96}.

The solution to this difficult problem was recently given by Maxim Kontsevich
\cite{Kt97} and this is in a way ``the cherry on the cake" of deformation
quantization (and contributed to getting him the Fields medal in 1998).
It involves very elaborate constructions, both conceptually
and computationally and makes an essential use of ideas coming from string
theory. We shall not attempt here to describe it in detail but give the
flavor of the development. The reader interested in more details should
refer to \cite{Kt97} and follow carefully subsequent developments. The
mathematically oriented reader may be interested in the ``Bourbaki-style"
presentation (by a French mathematician \cite{Oe98}, see quotation above...)
of the context and results.

Since, as we noted above, a Poisson bracket $P$ is a nontrivial 2-cocycle for
the Hochschild cohomology of the algebra $N=C^\infty(X)$, a natural question
is to decide whether this infinitesimal deformation can be extended to a
star product on $N$. Kontsevich answers positively, and more:
\begin{theorem}
Let $X$ be a differentiable manifold and $N=C^\infty(X)$. There is a
natural isomorphism between equivalence classes of deformations of
the null Poisson structure on $X$ and equivalence classes of differentiable
deformations of the associative algebra $N$; in particular, any Poisson
bracket $P$ on $X$ comes from a canonically defined (modulo equivalence)
star product.
\end{theorem}
With this concise formulation of the result (which gives a positive answer to
Kontsevich's formality conjecture) we see that, in this more general context,
a main result from the symplectic case is still valid: {\it classes
of star products correspond to classes of deformations of the Poisson
structure}. A deformation of the null Poisson structure is a formal
series $\Lambda(\hbar) = \sum_{n=1}^\infty \Lambda_n\hbar^n$ having vanishing
Schouten-Nijenhuis bracket with itself: $[\Lambda(\hbar),\Lambda(\hbar)]=0$.
Any given Poisson structure $\Lambda_0$ can be identified with the series
$\Lambda_0\hbar$.
\medskip

\noindent {\it Remark 3}. From Theorem 1 it is natural to conjecture that the
relation with \linebreak 1-differentiable deformations of the Poisson bracket
mentioned at the end of (II.2.1) extends to general Poisson manifolds, but
full proofs of that and a study of the relation with some 2-cohomology on the
manifold have not yet been given. Furthermore (and this is one of the
developments to come) a comparison of the proof with e.g. that of Fedosov
in the symplectic case is not done either. 
Finally, in his Berlin August 25, 1998 ``Fields medal" lecture, 
{\it Motivic Galois group and deformation quantization}, Kontsevich stressed
that the isomorphism described in Theorem~1 should be taken as one of a
family of isomorphisms and indicates that the motivic Galois group should
act on the moduli space of non-commutative algebras.
\medskip 

The bulk of the proof is the ``affine case", essentially when $X$ is some
${\Bbb R}^\ell$ (or an open set in it), the result being formulated in such
a way that ``gluing" charts, though still nontrivial, is not too difficult
to perform for an experienced mathematician. Doing so Kontsevich gives an
interesting explicit universal formula for the star product on such an $X$
where graphs (and Stokes formula) play a crucial r\^ole. The formula looks
like
\begin{equation}
u\ast v = \sum^{\infty}_{n=0}
\sum_{\Gamma\in G_n} w_\Gamma B_{\Gamma,\Lambda}(u,v)
\end{equation}
where $\Lambda$ is an arbitrary Poisson structure on an open domain in
${\Bbb R}^\ell$ and $G_n$ a set of labeled oriented graphs $\Gamma$.
The latter are pairs $(V_\Gamma,E_\Gamma)$ such that
$E_\Gamma \subset V_\Gamma\times V_\Gamma$ with $n+2$ vertices $V_\Gamma$ and
$2n$ edges $E_\Gamma$ satisfying some additional conditions (see \cite{Kt97},
section 2). $G_n$ has $(n(n+1))^n$ elements (1 for $n=0$).
$B_\Gamma$ is an explicitly defined bidifferential operator (of total order
$2n$, as in (4)) and $w_\Gamma$ a weight defined by an
absolutely convergent integral (of the exterior product of the differentials
of $2n$ harmonic angular variables associated with $\Gamma$, taken over the
space of configurations of $n$ numbered pairwise distinct points on the
Lobatchevsky upper half-plane).

\subsubsection*{II.2.3 Remarks on the infinite-dimensional case}
Poisson structures are known on infinite-dimensional manifolds since
a long time and there is an extensive literature on this subject, which alone
would require a book. A typical structure, for our purpose, is a symplectic
structure such as that defined by Segal \cite{Se74} (see also \cite{Ko74})
on the space of solutions of a classical field equation like
$\Box\Phi=F(\Phi)$, where $\Box$ is the d'Alembertian. Now if one considers
scalar valued functionals $\Psi$ over such a space of solutions, i.e. over
the phase space of initial conditions $\varphi(x)=\Phi(x,0)$ and
$\pi(x)={\frac{\partial}{\partial t}}\Phi(x,0)$, one can consider a Poisson
bracket defined by
\begin{equation} \label{PoiF}
P(\Psi_1,\Psi_2)=\int({\frac{\delta\Psi_1}{\delta\varphi}}
{\frac{\delta\Psi_2}{\delta\pi}}-{\frac{\delta\Psi_1}{\delta\pi}}
{\frac{\delta\Psi_2}{\delta\varphi}})dx
\end{equation}
where $\delta$ denotes the functional derivative. The problem is that while
it is possible to give a precise mathematical meaning to (\ref{PoiF}), the
formal extension to powers of $P$, needed to define the Moyal bracket,
is highly divergent, already for $P^2$.

The same difficulty is met if one takes e.g. a space $N$ of differentiable
functions on a Hilbert space with orthonormal basis
$\{p_k,q_k;k=1,\ldots,\infty\}$ and a Poisson bracket
$P(u,v)=\sum_{k=1}^\infty({\frac{\partial u}{\partial p_k}}
{\frac{\partial v}{\partial q_k}}-{\frac{\partial u}{\partial q_k}}
{\frac{\partial v}{\partial p_k}})$.
\smallskip

This is not so surprising for physicists who know from experience that
the correct approach to field theory is via normal ordering, and that
there are infinitely many inequivalent representations of the canonical
commutation relations (as opposed to the von Neumann uniqueness in the
finite-dimensional case, for projective representations).
Integral formulas, related to Feynman path integrals, can also be used with
some success. 
The participants at the \L \'od\'z meeting may be interested to learn that an
analogue of the pseudodifferential calculus in the infinite-dimensional case,
and especially the (``Wigner") notion of symbols of operators,
has been developed already in 1978 by Paul Kr\'ee and Ryszard R\c aczka 
\cite{KR78}.

We shall come back to this question in (II.3.2) with more specific
examples and give indications showing that with proper care the deformation
quantization approach can help making better mathematical sense of field
theory calculations done by theoretical physicists.

\subsubsection*{II.2.4 Generalized deformations, $n$-gebras and
related structures}
One of the mathematical reasons we started with the study of deformations
of Poisson brackets is related to the fact that it is the only one, among
classical infinite-dimensional algebras, which is not rigid, even at the
level of 1-differentiable deformations. In particular unimodular structures
(defined by a determinant) are rigid. It turns out that in connection with
Nambu mechanics, where the Poisson bracket is replaced by an $n$-bracket,
say a functional determinant, one meets structures of this type and it is
not a big surprise that a specific quantization (not of Heisenberg type)
was difficult to find.

Roughly speaking, a {\it generalized deformation} of a ${\Bbb{K}}$-algebra $A$
(associative, Lie or other) is a ${\Bbb{K}}$-algebra $A_\nu$ having $A$ for
limit as the deformation parameter $\nu\rightarrow 0$. Among the ``other"
algebras are of course the bialgebras to which we shall come back in (II.4.2)
when dealing with quantum groups (incidentally, since `al' means `the' in
Arabic, applied to a set containing only one element, the French denomination
`big\`ebre', imposed by Cartier, is far better).

Here we shall be concerned mostly with the so-called $n$-gebras, algebras $A$
endowed with a composition law  $A^n\rightarrow A$ satisfying some conditions
including skew-symmetry. Structures of this kind were introduced by Nambu
\cite{Nb73} in connection with his ``generalized mechanics" and
(in a paper published in an obscure journal) by Filippov \cite{Fi85}.
Serious interest in them developed only from 1992, when Takhtajan \cite{Ta94}
and independently Flato and Fr\o nsdal (unpublished) discovered that Nambu
$n$-brackets satisfy a generalization of Jacobi identity, called the
Fundamental Identity (FI); surprisingly enough, this identity had not been
discovered before. 

The resurgence of {\it operads} which occurred at the same time 
\cite{GK94,Lo94} and are related to $n$-gebras \cite{Gn95}, as well
as the new notion of {\it strong homotopy Lie algebras} introduced then by 
Stasheff and is also related to deformation theory \cite{St97} add to the
interest in these structures.

Recently, there have been several works dealing with various generalizations
of Poisson structures by extending the binary bracket to an $n$-bracket.
The main point for these generalizations is to look for the corresponding
identity which would play the r\^ole of Jacobi identity for the usual Poisson
bracket. Indeed, in view of generalizations, the Jacobi identity for a Lie
$2$-bracket can be presented in a number of ways \cite{Ga96} among which two
have been recently extensively studied. The most straightforward way is to
require that the sum over the symmetric group ${\frak{S}}_3$ of the composed
brackets $[[\cdot,\cdot],\cdot]$ is zero. When extended to $n$-brackets,
leads to the notion of {\it generalized Poisson structures} studied e.g. in
\cite{APP96}; the corresponding identity is obtained by complete
skew-symmetrization of the $2n-1$ composed brackets when $n$ is even;
this is equivalent to require that the Schouten bracket of the
$n$-tensor defining the $n$-bracket with itself vanishes.

A physically more appealing way is to say that the adjoint map
$b\mapsto [a,b]$ is a Lie algebra derivation. Indeed this means that the
bracket of conserved quantities is again a conserved quantity. The two
formulations coincide only for $n=2$ and for $n\geq3$ the latter is a
stronger requirement. This {\it Fundamental Identity} of Nambu Mechanics can
be written:
\begin{eqnarray}\label{FI}
&& [x_1,\ldots,x_{n-1}, [y_1,\ldots,y_{n}]]\cr
&&\quad - \sum_{1\leq i \leq n}[y_1,\ldots,y_{i-1},y_{i+1},\ldots,y_{n},
[x_1, \ldots,x_{i-1},y_i,x_{i},\dots,x_{n-1}]]=0.
\end{eqnarray}
{\it Nambu brackets} (like Poisson brackets and commutators such as the Moyal
bracket, for $n=2$) are $n$-brackets required to satisfy, with respect to the
usual algebra multiplication and in addition to skew-symmetry
$\{x_1,\ldots,x_n\} = \epsilon(\sigma) \{x_{\sigma_1},\ldots,x_{\sigma_n}\}$
$\forall \sigma\in {\frak{S}}_n$
and the FI, a {\it Leibniz rule}:
\begin{equation}\label{LR}
\{x_0x_1,x_2,\ldots,x_n\}=x_0\{x_1,x_2,\ldots,x_n\}+\{x_0,x_2,\ldots,x_n\}x_1.
\end{equation}
The related cohomologies are not yet completely known, though a major step
in this direction was done in \cite{Ga96} where one can also find a very
interesting and detailed study of all intermediate possibilities
between the two generalizations described here (generalized Poisson and
Nambu). In (II.3.3) we shall nevertheless indicate, specializing to the case
$A=N$, how one can quantize Nambu brackets using generalized deformations
based on the factorization of polynomials and methods of second quantization.
One of the steps there (and this produces a non-DrG-deformation) is an
operation, the effect of which is that in products the deformation parameter
$\hbar$ behaves as if it was nilpotent (e.g. multiplied by a Dirac $\gamma$
matrix).

This last fact has very recently induced Pinczon \cite{Pi97} and Nadaud
\cite{Na98} to generalize Gerstenhaber theory to the case of a deformation
parameter $\sigma$ which {\it does not commute with the algebra}.
A similar theory can be done in this case, with appropriate
cohomologies. While that theory does not reproduce the above mentioned
Nambu quantization, it gives new and interesting results.
For instance \cite{Pi97}, while the Weyl algebra $W_1$
(generated by the Heisenberg Lie algebra ${\frak{h}}_1$) is known \cite{DC85}
to be DrG-rigid, it can be nontrivially deformed in such a {\it supersymmetric
deformation theory} to the supersymmetry enveloping algebra
${\cal{U}}({\frak{osp}}(1,2))$; or \cite{Na98}, on the polynomial algebra
${\Bbb C}[x,y]$ in 2 variables, Moyal-like products of a new type were
discovered. This is another example of a motivated study which goes beyond
a generally accepted framework.

\subsection*{II.3 Physical Applications}
In this subsection and the following, I shall present a few of the numerous
developments which have made use of deformation quantization and/or are
strongly related to it. The presentation made, and therefore the
bibliography, is by no means exhaustive -- more than a whole volume would
be needed for that -- and the absence of reference to any specific work
does not (in general) reflect a lack of appreciation; I did not even quote
all of my publications in the domain. The aim of these two last subsections
(in fact, of all this review) is mainly to give the flavor of the many
facets of deformation quantization and quite naturally the presentation
will be somewhat biased towards the works of our group. Nevertheless
the interested reader should be able to complete whatever is missing by
a kind of ``hyper-referencing", looking at references of references a few
times.

\subsubsection*{II.3.1 Quantum mechanics}
Let us start with a phase space $X$, a symplectic (or Poisson) manifold and
$N$ an algebra of classical observables (functions, possibly including
distributions if proper care is taken for the product). We shall call
{\it star quantization} a star product on $N$ invariant (or sometimes only
covariant) under some Lie algebra ${\frak{g}}_0$ of ``preferred observables".
Invariance of the star product ensures that the classical and quantum
evolutions of observables under a Hamiltonian $H \in {\frak{g}}_0$ will
coincide \cite{BFFLS78}. The typical example is the Moyal product on
$W = {\Bbb{R}}^{2\ell}$.
\medskip

\noindent {\it II.3.1.1 Spectrality.} Physicists want to get numbers
matching experimental results, e.g. for energy levels of a system. That is
usually achieved by describing the spectrum of a given Hamiltonian $\hat{H}$
supposed to be a self-adjoint operator so as to get a real spectrum and
so that the evolution operator (the exponential of $it \hat{H}$)
is unitary (thus preserves probability). A similar spectral theory can
be done here, in an {\it autonomous manner}. The most efficient way to
achieve it is to consider \cite{BFFLS78} the {\it star exponential}
(corresponding to the evolution operator)
\begin{equation}\label{starexp}
\mbox{Exp}(Ht) \equiv \sum^\infty_{n=0}
\frac{1}{n!}\left(\frac{t}{i\hbar}\right)^n(H*)^n
\end{equation}
where $(H*)^n$ means the $n^{\mbox{\footnotesize th}}$ star power of the
Hamiltonian $H \in N$ (or $N[[\nu]]$).
Then one writes its Fourier-Stieltjes transform $d\mu$ (in the distribution
sense) as $\mbox{Exp}(Ht) = \int e^{\lambda t/i\hbar} d\mu (\lambda)$ and
defines {\it the spectrum of $(H /\hbar)$} as the {\it support $S$ of $d\mu$}
(incidentally this is the definition given by L. Schwartz for the spectrum
of a distribution, out of motivations coming from Fourier analysis).
In the particular case when $H$ has discrete spectrum, the integral can
be written as a sum (see the top equation in (\ref{spec}) below for a
typical example): the distribution $d\mu$ is a sum of ``delta functions"
supported at the points of $S$ multiplied by the symbols of the
corresponding eigenprojectors.

In different orderings with various weight functions $w$ in (\ref{weyl}) one
gets in general different operators for the same classical observable $H$,
thus different spectra. For $X = {\Bbb{R}}^{2\ell}$ all orderings are
mathematically equivalent (to Moyal under the Fourier transform $T_w$
of the weight function $w$). This means that every observable $H$ will have
the same spectrum under Moyal ordering as $T_w H$ under the equivalent
ordering. But this does not imply physical equivalence, i.e. the fact that
$H$ will have the same spectrum under both orderings. In fact the opposite
is true \cite{CFGS85}: if two equivalent star products are isospectral
(give the same spectrum for a large family of observables and all $\hbar$),
they are identical.

It is worth mentioning that our definition of spectrum permits to define
a spectrum even for symbols of non-spectrable operators, such as the
derivative on a half-line which has different deficiency indices; this
corresponds to an infinite potential barrier (see also \cite{UU96}
for detailed studies of similar questions). That is one of the many
advantages of our autonomous approach to quantization.
\medskip

\noindent {\it II.3.1.2 Applications.}  In quantum mechanics it is preferable
to work (for $X = {\Bbb{R}}^{2\ell}$) with the star product that has maximal
symmetry, i.e. ${\frak{sp}}({\Bbb{R}}^{2\ell})\cdot{\frak{h}}_\ell$
as algebra of preferred observables: the Moyal product. One indeed
finds \cite{BFFLS78} that the star exponential of these observables
(polynomials of order $\leq 2$) is proportional to the usual exponential.
More precisely, if $H=\alpha p^2+\beta pq + \gamma q^2 \in {\frak{sl}}(2)$
with $p,q \in  {\Bbb{R}}^{\ell}$, $\alpha,\beta, \gamma \in {\Bbb{R}}$,
setting $d = \alpha\gamma -\beta^2$ and $\delta = \vert d \vert^{1/2}$
one gets (the sums and integrals appearing in the various expressions of
the star exponential being convergent as distributions, both in phase-space
variables and in $t$ or $\lambda$)
\begin{equation}\label{expalt}
\mbox{Exp}(Ht) = \left\{ \begin{array}{r@{\quad \mbox{for} \quad}l}
 (\cos\delta t)^{-l} \exp((H/i\hbar\delta)\tan(\delta t))& d > 0 \\
\exp(Ht/i\hbar)& d = 0 \\
(\cosh\delta t)^{-l}\exp((H/i\hbar\delta)\tanh(\delta t))& d < 0
\end{array}\right.
\end{equation}
hence the Fourier decompositions
\begin{equation}\label{spec}
\mbox{Exp}(Ht) = \left\{\begin{array}{r@{\quad \mbox{for} \quad}l}
\sum^\infty_{n=0} \Pi^{(\ell)}_n e^{(n + {\ell \over 2})t} & d > 0 \\
\int^\infty_{- \infty} e^{\lambda t/i\hbar} \Pi(\lambda, H)d\lambda & d < 0
\end{array}\right.
\end{equation}
We thus get the discrete spectrum $(n + {\ell\over 2})\hbar$ of the
{\it harmonic oscillator} and the continuous spectrum ${\Bbb{R}}$
for the dilation generator $pq$. The eigenprojectors $\Pi^{(\ell)}_n$ and
$\Pi (\lambda, H)$ are given \cite{BFFLS78} by known special functions on
phase-space (generalized Laguerre and hypergeometric, multiplied by some
exponential). Formulas (\ref{expalt}) and (\ref{spec}) can, by analytic
continuation, be given a sense outside singularities and even (as
distributions) for singular values of $t$.

Other examples can be brought to this case by functional manipulations
\cite{BFFLS78}. For instance the Casimir element $C$ of ${\frak{so}}(\ell)$
representing {\it angular momentum}, which can be written
$C = p^2 q^2 - (pq)^2 - \ell(\ell - 1){\hbar^2 \over 4}$,
has $n (n + (\ell-2)) \hbar^2$ for spectrum. For the {\it hydrogen atom},
with Hamiltonian $H = {1 \over 2} p^2 - \vert q \vert^{-1}$, the Moyal
product on ${\Bbb{R}}^{2\ell+2}$ ($\ell=3$ in the
physical case) induces a star product on $X = T^*S^\ell$; the energy
levels, solutions of $(H-E) * \phi = 0,$ are found from (\ref{spec})
and the preceding calculations for angular momentum to be (as they should,
with $\ell=3$) $E = {1 \over 2} (n+1)^{-2} \hbar^{-2}$ for the discrete
spectrum, and $E \in {\Bbb{R}}^+$ for the continuous spectrum.

We thus have recovered, in a completely autonomous manner entirely within
deformation quantization, the results of ``conventional" quantum mechanics
in those typical examples (and many more can be treated similarly).
It is worth noting that the term ${\ell \over 2}$ in the harmonic oscillator
spectrum, obvious source of divergences in the infinite-dimensional case,
disappears if the normal star product is used instead of Moyal -- which
is one of the reasons it is preferred in field theory.
\medskip

\noindent{\it II.3.1.3 Remark on convergence}. We have always considered
star products as formal series and looked for convergence only in specific
examples, and then generally in the sense of distributions. The same applies
to star exponentials, as long as each coefficient in the formal series is
well defined. In the case of the harmonic oscillator or more generally
the preferred observables $H$ in Weyl ordering, this study was
facilitated by the fact that the powers $(H*)^n$ are polynomials in $H$.
Moreover, in the case of star exponentials, a notion of convergence stronger
than as distributions would require considerations analogous to the problem
of analytic vectors in Lie groups representations \cite{FSSS72} and pose
problems also when looking at their Fourier decomposition. Nevertheless
some authors (see e.g. \cite{Ri97} and references therein) insist in making
a stronger parallel with operator algebras and look for domains (in $N$)
where one has convergence, speaking of ``strict" deformation quantization.
In the Weyl case on ${\Bbb R}^{2\ell}$ this question is related to the
``Wigner image" of classes of operators \cite{Ha84,Ka86,Ma86} on
$L^2({\Bbb R}^{\ell})$. While this is a perfectly legitimate mathematical
problem, we do not feel that it is physically wise. In particular it lacks
the flexibility that exists in deformation quantization as compared with
physics based on algebras of bounded operators or on the less developed
approach of algebras of unbounded operators \cite{Ju85}, and puts deformation
quantization back into the $C^*$~algebras Procrustean bed. Therefore
the question of ``strict" convergence shall not be asked, except in examples.
This does not mean that one should not look for domains where pointwise
convergence can be proved; this was done e.g. for Hermitian symmetric spaces
\cite{CGR95}. But it should be clearly understood that one can consider
wider classes of observables -- in fact, the latter tend to be physically more
interesting. It is also worthwhile to adapt to deformation quantization
procedures that were successful in operator algebras, like the GNS
construction \cite{BW96,BNW97}.

\subsubsection*{II.3.2 Path integrals, field theory and statistical mechanics}
\noindent{\it II.3.2.1 Path integrals} are intimately connected to star
exponentials. In fact, in quantum mechanics the path integral of the action
is nothing but the partial Fourier transform of the star exponential
(\ref{starexp}) with respect to the momentum variables, for
$X = {\Bbb{R}}^{2\ell}$ as phase space with the Moyal star product
\cite{PS79}. For normal ordering the path integral is essentially the star
exponential \cite{Di90} and we shall come back to it in (II.3.2.2).

For compact groups the star exponential $E$ (defined in a similar manner,
see below (II.4.1)) can be expressed in terms of unitary characters using
a global coherent state formalism \cite{CN91} based on the Berezin
dequantization of compact group representation theory used in
\cite{ACG88,Mr86} (it gives star products somewhat similar to normal
ordering); the star exponential of any Hamiltonian on $G/T$ (where $T$
is a maximal torus in the compact group $G$) is then equal to the path
integral for this Hamiltonian.
\medskip

\noindent{\it II.3.2.2 Field theory}.
The deformation quantization of a given classical field theory consists in
the giving a proper definition for a star product on the infinite-dimensional
manifold of initial data for the classical field equation (see II.2.3) and
constructing with it, as rigorously as possible, whatever physical
expressions are needed.
As in other approaches to field theory, here also one faces serious divergence
difficulties as soon as one is considering interacting fields theory, and
even at the free field level if one wants a mathematically rigorous theory.
But the philosophy in dealing with the divergences is significantly different
and one is in position to take advantage of the cohomological features of
deformation theory to perform what can be called {\it cohomological
renormalization}.

Starting with some star product $\ast$ (e.g. an infinite-dimensional version
of a Moyal-type product or, better, a star product similar to the
normal star product (\ref{np})) on the manifold of initial data, one would
interpret various divergences appearing in the theory in terms of
coboundaries (or cocycles) for the relevant Hochschild cohomology.
Suppose that we are suspecting that a term in a cochain of the product $\ast$
is responsible for the appearance of divergences. Applying the procedure
described in (I.4.2.2), we can try to eliminate it, or at least get a lesser
divergence, by subtracting at the relevant order a coboundary; we would then
get a better theory with a new star product, equivalent to the original one.
Furthermore, since in this case we can expect to have at each order an
infinity of non equivalent star products, we can try to subtract a cocycle
and then pass to a non equivalent star product whose lower order cochains
are identical to those of the original one. We would then make an analysis of
the divergences up to order $\hbar^r$, identify a divergent cocycle, remove
it, and continue the procedure (at the same or hopefully a higher order).
Along the way one should preserve the usual properties of a quantum field
theory (Poincar\'e covariance, locality, etc.) and the construction of
adapted star products should be done accordingly. The complete implementation
of this program should lead to a cohomological approach to renormalization
theory. 

A very good test for this approach would be to start from classical
electrodynamics, where (among others) the existence of global
solutions and a study of infrared divergencies were recently rigorously
performed \cite{FST97}, and go towards mathematically rigorous QED.
Physicists will think that spending so much effort in trying to give complete 
mathematical sense to recipes that work so well is a waste of time, but
I am sure that the mathematical tools needed will prove very efficient.
Would De~Gaulle have been a mathematician he could have said about
this scheme ``vaste programme" 
(supposedly his answer to a minister
who wanted to get rid of all stupid bureaucrats); but had he been a scientist
he would probably have been a physicist and share the attitude of too many
physicists towards mathematics: ``l'intendance suivra" 
(``Supply Corps will follow", needed logistics will be provided).
 
In the case of free fields, one can write down an explicit expression for
a star product corresponding to normal ordering. Consider a (classical)
free massive scalar field $\Phi$ with initial data $(\phi,\pi)$ in the
Schwartz space ${\cal{S}}$. The initial data  $(\phi,\pi)$ can
advantageously be replaced by their Fourier modes $({\bar{a}},a)$ which after
quantization become the usual creation and annihilation operators,
respectively. The normal star product $\ast_N$ is formally equivalent
to the Moyal product and an integral representation for $\ast_N$ is given by:
\begin{equation} \label{np}
(F\ast_N G)(\bar{a},a)=\int_{{\cal{S}}'\times {\cal{S}}'}
d\mu(\bar{\xi},\xi) F(\bar{a},a +\xi) G(\bar{a} + \bar{\xi}, a),
\end{equation}
where $\mu$ is a Gaussian measure on ${\cal{S}}'\times {\cal{S}}'$ defined by
the characteristic function $\exp(-\frac{1}{\hbar}\int dk\ {\bar{a}}(k)a(k))$
and $F,G$ are holomorphic functions with semi-regular kernels. Likewise,
Fermionic fields can be cast in that framework by considering functions
valued in some Grassmann algebra and super-Poisson brackets (for the
deformation quantization of the latter see e.g. \cite{Bm96}).

For the normal product (\ref{np}) one can formally consider interacting fields.
It turns out that the star exponential of the Hamiltonian is, up to a
multiplicative well-defined function, equal to Feynman's path integral.
For free fields, we have a mathematical meaningful equality between
the star exponential and the path integrals as both of them are defined
by a Gaussian measure, and hence well-defined. In the interacting fields
case, giving a rigorous meaning to either of them would give a meaning
to the other.

The interested reader will find in \cite{Di90} calculations performing some
steps in the above direction, for free scalar fields and the Klein-Gordon
equation, and an example of cancellation of some infinities in
$\lambda\phi_2^4$-theory via a $\lambda$-dependent star product equivalent
to a normal star product. Finally we mention here for 
more completeness (though neither is directly related to what precedes) the
symbolic calculus of \cite{KR78} and the Fedosov-like
approach to self-dual Yang-Mills and gravity of \cite{GPP98}.
\medskip

\noindent{\it II.3.2.3 Statistical mechanics}. In view of our philosophy on
deformations, a natural question to ask is their {\it stability}: Can
deformations be further deformed, or does ``the buck stops there"? As we
indicated at the beginning of this review and shall exemplify with quantum
groups, the answer to that question may depend on the context. Here is
another example.

If one looks for deformations of the Poisson bracket Lie algebra $(N,P)$
one finds (assuming mild technical assumptions on parity of cochains in
\cite{BL81} which, in view of the classification of star products, are
not required) that a further deformation of the Moyal bracket, with another
deformation parameter $\rho$, is again a Moyal bracket for a $\rho$-deformed
Poisson structure; in particular, for $X={\Bbb R}^{2\ell}$,
{\it quantum mechanics viewed as a deformation is unique and stable}.

Now, for the associative algebra $N$, the only {\it local} associative
composition law is \cite{Ru84} of the form $(u,v)\mapsto ufv$ for some
$f\in N$. If we take $f=f_\beta\in N[[\beta]]$ we get a 0-differentiable
deformation (with parameter $\beta$) of the usual product, which for
convenience we shall call here a Rubio product. We were thus lead
\cite{BFLS84} to look, starting from a $\ast_\nu$ product, for a new
composition law
\begin{equation}\label{tilstar}
(u,v)\mapsto u{\tilde{\ast}}_{\nu,\beta}v=u\ast_\nu f_{\nu,\beta}\ast_\nu v
\quad \mbox{with} \quad f_{\nu,\beta}=\sum_{r=0}^\infty \nu^{2r}f_{2r,\beta}
\in N[[N[[\nu^2]],\beta]]
\end{equation}
where $f_{0,\beta}\equiv f_\beta \neq 0$ and $f_0=1$. The transformation
$u\mapsto T_{\nu,\beta}u=f_{\nu,\beta}\ast_\nu u$ intertwines $\ast_\nu$
and ${\tilde{\ast}}_{\nu,\beta}$ but it is not an equivalence of star products
because ${\tilde{\ast}}_{\nu,\beta}$ is not a star product: it is a
($\nu,\beta$)-deformation of the usual product (or a $\nu$-deformation of the
Rubio product) with at first order in $\nu$ the driver
given by $P_\beta(u,v)=f_\beta P(u,v)+uP(f_\beta,v)-P(f_\beta,u)v$, a
conformal Poisson bracket associated with a {\it conformal symplectic
structure} given by the 2-tensor $\Lambda_\beta=f_\beta\Lambda$ and the
vector $E_\beta=[\Lambda,f_\beta]$.

In view of applications we suppose given a star product, denoted $\ast$, on
some algebra ${\cal{A}}$ of observables (possibly defined on some
infinite-dimensional phase-space) and take for $f_{\nu,\beta}$ the
exponential $g_\beta\equiv \exp_\ast(c\beta H)=
1+\sum_{n=1}^\infty \frac{(c\beta)^n}{n!}(H\ast)^n$ with $c=-{\frac{1}{2}}$
(we omit $\nu$ from now on and write ${\tilde{\ast}}$ for
${\tilde{\ast}}_{\nu,\beta}$). The star exponential $\mbox{Exp}(Ht)$ defines
an automorphism $u\mapsto\alpha_t(u)=\mbox{Exp}(-Ht)\ast u\ast\mbox{Exp}(Ht)$.
A {\it KMS state} $\sigma$ on ${\cal{A}}$ is a state (linear functional)
satisfying, $\forall a,b\in {\cal{A}}$, the Kubo--Martin--Schwinger condition
$\sigma(\alpha_t(a)\ast b)=\sigma(b\ast\alpha_{t+i\hbar\beta}(a))$.
Then the (quantum) KMS condition can be written \cite{BFLS84}, with
$[a,b]_\beta=(i/\hbar)(a{\tilde{\ast}}b-b{\tilde{\ast}}a)$, simply
$\sigma(g_{-\beta}\ast[a,b]_\beta)=0$: up to a conformal factor, a KMS state
is like a trace with respect to this new product. The (static) classical
KMS condition is the limit for $\hbar=0$ of the quantum one. So we can
recover known features of statistical mechanics by introducing a new
deformation parameter $\beta=(kT)^{-1}$ and the related conformal symplectic
structure. This procedure commutes with usual deformation quantization.
Finally let us mention that recently several people \cite{BRW98,Ws97} have
considered the question of KMS states and related modular automorphisms
from a more conventional point of view in deformation quantization.

\subsubsection*{II.3.3 Nambu mechanics and its quantization}
We mention this aspect here mainly for the sake of completeness, as an example
of generalized deformation. A somewhat detailed recent review can be found in
\cite{FDS97} (see also \cite{Fl98}), so we shall just briefly indicate a few
highlights.

Nambu \cite{Nb73} started with a kind of ``Hamilton equations" on ${\Bbb R}^3$,
of the form $\frac{d{r}}{dt}=\nabla g({r}) \wedge \nabla h({r})$,
$r=(x,y,z)\in {\Bbb R}^3$, where $x$, $y$, $z$ are the dynamical variables
and $g$, $h$ are two functions of ${r}$. Liouville theorem follows directly
from the identity
$\nabla\cdot(\nabla g({r}) \wedge \nabla h({r}))=0$,
which tells us that the velocity field in the above equation is 
divergenceless. From this we derive the evolution of a function $f$ on
${\Bbb R}^3$:
\begin{equation}\label{Jac}
\frac{df}{dt}=\frac{\partial(f,g,h)}{\partial(x,y,z)} ,
\end{equation}
where the right-hand side is the Jacobian of the mapping
${\Bbb R}^3 \rightarrow {\Bbb R}^3$ given by $(x,y,z)\mapsto (f,g,h)$.
In this ``baby model for integrable systems", Euler equations
for the angular momentum of a rigid body are obtained when the dynamical
variables are taken to be the components of the angular  momentum vector
$ L = (L_x, L_y, L_z)$, $g$ is the total kinetic energy
${L_x^2 \over 2I_x} + {L_y^2 \over 2I_y} + {L_z^2 \over 2I_z}$
and $h$ the square of the angular momentum $L_x^2 + L_y^2 + L_z^2$.
Other examples can be given, in particular Nahm's equations for static
${\frak {su}}(2)$ monopoles, $\dot{x}_i = x_j x_k$ ($i,j,k = 1,2,3$) in
${\frak {su}}(2)^* \sim {\Bbb R}^3$, with $h = x_1^2 -x_2^2$,
$g = x_1^2 -x_3^2$, etc.
Here the principle of least action, which states that the classical
trajectory $C_1$ is an extremal of the action functional
$ A(C_1) = \int_{C_1} (pdq - Hdt) $, is replaced  by a similar one
\cite{Ta94} with a 2-dimensional cycle $C_2$ and ``action functional''
 $A(C_2) =  \int_{C_2} (xdy \wedge dz - h dg \wedge dt) $
(which bears some flavor of strings and some similitude with
the cyclic cocycles of Connes \cite{Co94}).

Expression (\ref{Jac}) was easily generalized to $n$ functions $f_i$,
$i=1, \ldots, n$. One introduces an $n$-tuple of functions on ${\Bbb R}^n$
with composition law given by their Jacobian, linear canonical
transformations ${\rm{SL}}(n,{\Bbb R})$ and a corresponding $(n-1)$-form
which is the analogue of the Poincar\'e-Cartan integral invariant.
The Jacobian has to be interpreted as a generalized Poisson bracket:
It is skew-symmetric with respect to the $f_i$'s, satisfies the FI which is
an analogue of the Jacobi identity (but was discovered much after \cite{Nb73})
and a derivation of the algebra of smooth functions on  ${\Bbb R}^n$
(i.e., the Leibniz rule is verified in each argument, e.g.
$\{f_1f_2, f_3,\ldots, f_{n+1}\}=f_1\{f_2,\ldots,f_{n+1}\}+
 \{f_1, f_3, \ldots, f_{n+1}\}f_2$, etc.).
Hence there is a complete analogy with the Poisson bracket formulation
of Hamilton equations, including the important fact that the components of the
$(n-1)$-tuple of ``Hamiltonians" $(f_2,\ldots ,f_n)$ are constants of motion.

Shortly afterwards it was shown \cite{BF75,MS76} that Nambu mechanics could
be seen as a coming from constrained Hamiltonian mechanics; e.g. for
${\Bbb R}^3$ one starts with ${\Bbb R}^6$ and an identically vanishing
Hamiltonian, takes a pair of second class constraints to reduce it to some
${\Bbb R}^4$ and one more first-class Dirac constraint, together with time
rescaling, will give the reduction. This ``chilled" the domain for almost
20 years -- and gives a physical explanation to the fact that Nambu could
not go beyond Heisenberg quantization.

In order to quantize the Nambu bracket, a natural idea is to replace, in the
definition of the Jacobian, the pointwise product of functions by a deformed
product. For this to make sense, the deformed product should be
Abelian, so we are lead to consider commutative DrG-deformations of an
associative and commutative product. Looking first at polynomials
(this restriction can be removed \cite{PG97}) we are lead to the commutative
part of Hochschild cohomology called Harrison cohomology, which is trivial
\cite{Ba68,GS88}.
Dealing with polynomials, a natural idea is to factorize them and take
symmetrized star products of the factors. More precisely we introduce
an operation $\alpha$ which maps a product of factors into a symmetrized
tensor product (in a kind of Fock space) and an evaluation map $T$
which replaces tensor product by star product. Associativity will be
satisfied if $\alpha$ annihilates the deformation parameter $\hbar$
(there are still $\hbar$-dependent terms in a product due to the last action
of $T$); intuitively one can think of a deformation parameter which is
$\hbar$ times a Dirac $\gamma$ matrix. This fact brought us to generalized
deformations, but even this was not enough. Dealing with distributivity of
the product with respect to addition and with derivatives posed difficult
problems. In the end we took for observables Taylor developments of elements
of the algebra of the semi-group generated by irreducible polynomials
(``polynomials over polynomials", inspired by second quantization techniques)
and were then able to perform a meaningful quantization of these Nambu-Poisson
brackets (cf. \cite{DFST97} for more details and \cite{DF97} for subsequent
development).

\subsection*{II.4 Related Mathematical Developments}
\subsubsection*{II.4.1 Star representation theory of Lie groups}
Let $G$ be a Lie group (connected and simply connected), acting by
symplectomorphisms on a symplectic manifold $X$ (e.g. coadjoint orbits
in the dual of the Lie algebra ${\frak{g}}$ of $G$). The elements
$x,y \in {\frak{g}}$ will be supposed realized by functions $u_x, u_y$ in $N$
so that their Lie bracket $[x,y]_{\frak g}$ is realized by $P(u_x, u_y)$.
Now take a $G$-covariant star-product $*$, that is
$P(u_x,u_y) = [u_x, u_y]\equiv (u*v-v*u)/2\nu$, which shows that the map
${\frak{g}}\ni x \mapsto (2\nu)^{-1}u_x\in N$ is a Lie algebra morphism.
The appearance of $\nu^{-1}$ here and in the trace (see (I.1)) cannot be
avoided and explains why we have often to take into account both $\nu$ and
$\nu^{-1}$. We can now define the {\it star exponential}
\begin{equation}\label{stargp}
E (e^x) = {\mbox{Exp}}(x) \equiv\sum^\infty_{n=0} (n!)^{-1} (u_x /2 \nu)^{*n}
\end{equation}
where $x \in {\frak{g}}$, $e^x \in G$ and the power $*n$ denotes the
$n^{\mbox{\footnotesize th}}$ star-power of the corresponding function. By the
Campbell-Hausdorff formula one can extend $E$ to a {\it group homomorphism}
$E : G \rightarrow (N[[\nu,\nu^{-1}]], *)$ where, in the formal series,
$\nu$ and $\nu^{-1}$ are treated as independent parameters for the time being.
Alternatively, the values of $E$ can be taken in the algebra
$({\cal P} [[\nu^{-1}]], *)$, where ${\cal P}$ is the algebra generated by
${\frak g}$ with the $*$-product (a representation of the enveloping algebra).

We call {\it star representation} \cite{BFFLS78,Fr78} of $G$ a distribution
${\cal E}$ (valued in ${\mbox{Im}} E$) on $X$ defined by
$D \ni f \mapsto {\cal E} (f) = \int_G f(g) E(g^{-1}) dg$
where $D$ is some space of test-functions on $G$. The corresponding
{\it character} $\chi$ is the (scalar-valued) distribution defined by
$D \ni f \mapsto \chi (f) = \int_X {\cal E} (f)d\mu $, $d\mu$ being a
quasi-invariant measure on $X$.

The character is one of the tools which permit a comparison with usual
representation theory. For semi-simple groups it is singular at the origin
in irreducible representations, which may require caution in computing the
star exponential (\ref{stargp}). In the case of the harmonic oscillator that
difficulty was masked by the fact that the corresponding representation of
${\frak{sl}}(2)$ generated by $(p^2,q^2,pq)$ is integrable to a double
covering of ${\rm{SL}}(2,{\Bbb{R}})$ and decomposes into a sum
$D({\frac{1}{4}})\oplus D({\frac{3}{4}})$: the singularities at the origin
cancel each other for the two components.

This theory is now very developed, and parallels in many ways the usual
(operatorial) representation theory. It is not possible here to give
a detailed account of all of them, but among notable results one may quote:
\smallskip

\noindent i) 
An exhaustive treatment of {\it nilpotent} or solvable exponential
\cite{AC85} and even {\it general solvable} Lie groups \cite{ACL95}.
The coadjoint orbits are there symplectomorphic to ${\Bbb R}^{2\ell}$
and one can lift the Moyal product to the orbits in a way that is adapted
to the Plancherel formula. Polarizations are not required, and
``star-polarizations" can always be introduced to compare with usual theory.
Wavelets \cite{Da95}, important in signal analysis, are manifestations of
star products on the (2-dimensional solvable) affine group of ${\Bbb{R}}$
or on a similar 3-dimensional solvable group \cite{BB98}.
\smallskip

\noindent ii) 
For {\it semi-simple} Lie groups an array of results is already
available, including \cite{ACG88,Mr86} a complete treatment of the
{\it holomorphic discrete series} (this includes the case of compact
Lie groups) using a kind of Berezin dequantization, and scattered results
for specific examples. Similar techniques have also been used 
\cite{CGR95,Ka98} to find invariant star products on K\"ahler and 
Hermitian symmetric spaces (convergent for an appropriate dense subalgebra). 
Note however, as shown by recent developments of unitary representations
theory (see e.g. \cite{Sch97}),
that for semi-simple groups the coadjoint orbits alone are no more sufficient
for the unitary dual and one needs far more elaborate constructions.
\smallskip

\noindent iii)
For semi-direct products, and in particular the Poincar\'e and
Euclidean groups, an autonomous theory has also been developed (see e.g.
\cite{ACM81}).
\smallskip

Comparison with the usual results of ``operatorial" theory of Lie group
representations can be performed in several ways, in particular by
constructing an invariant Weyl transform generalizing (\ref{weyl}),
finding ``star-polarizations" that always exist, in contradistinction
with the geometric quantization approach (where at best one can
find complex polarizations), study of spectra (of elements in the center
of the enveloping algebra and of compact generators) in the sense of
(II.3.1), comparison of characters, etc. Note also in this context that the
pseudodifferential analysis and (non autonomous) connection with
quantization developed extensively by Unterberger, first in
the case of ${\Bbb{R}}^{2\ell}$, has been recently extended to the above
invariant context \cite{UU96}. But our main insistence is that
the theory of star representations is an {\it autonomous} one that can be
formulated completely within this framework, based on coadjoint orbits
(and some additional ingredients when required).

\subsubsection*{II.4.2 Quantum groups}
Around 1980 Kulish and Reshetikhin \cite{KR81}, for purposes related
to inverse scattering and 2-dimensional models, discovered a strange
modification of the ${\frak{sl}}(2)$ Lie algebra, where the commutation
relation of the two nilpotent generators is a sine in the semi-simple
generator instead of being a multiple of it -- this in fact requires some
completion of the enveloping algebra ${\cal{U}}({\frak{g}})$. The theory was
developed in the first half of the 80's by the Leningrad school of L. Faddeev
\cite{FRT90}, systematized by V. Drinfeld who developed the Hopf algebraic
context and coined the extremely effective (though somewhat misleading) term
of {\it quantum group} \cite{Dr86} and from the enveloping algebra point of
view by Jimbo \cite{Ji85}. Shortly afterwards, Woronowicz \cite{Wn87}
realized these models in the context of the noncommutative geometry of
Alain Connes \cite{Co94} by matrix pseudogroups, with coefficients
(satisfying some relations) in $C^*$ algebras.
A typical example of such Hopf algebras is a Poisson Lie group, a Lie group
$G$ with compatible Poisson structure i.e. a Poisson bracket $P$ on
$N = C^{\infty}(G)$, considered as a bialgebra with coproduct defined by
$\Delta u(g,g')=u(gg')$, $g,g' \in G$, satisfying
$\Delta P (u,v) = P(\Delta u,\Delta v)$, $ u,v \in N$.

Now the topological dual of $N$ is the space $N'$ of distributions with
compact support on $G$; it includes $G$ (Dirac's $\delta$s at the points of
$G$) and a completion of ${\cal{U}}({\frak{g}})$ (differential operators).
Taking an adequate subspace $N_0$ of $N$ (generated by the coefficients of
suitably chosen representations, e.g. the ``well-behaved" vectors of
Harish Chandra) will give a dual $N_0'\supset N'$. All these are reflexive
(the bidual coincides with the original space; the algebraic dual of a
Hopf algebra is in general not a Hopf algebra). This is the basis of the
theory of topological Hopf algebras developed recently, first for $G$ compact
\cite{BFP92,BFGP94} and then for $G$ semi-simple and in general \cite{BP95}.
In the compact or semi-simple case the quantum group is obtained by giving
a star product on $N$ or $N_0$ and keeping unchanged the coproduct (what is
called a preferred deformation) or equivalently by deforming the coproduct
in the dual (and keeping the product unchanged). Associativity of the star
product corresponds to the Yang-Baxter equation, and the Faddeev-Reshetikhin-%
Takhtajan and Jimbo models of quantum groups can be seen in this way.
Also, all Poisson-Lie groups can be quantized \cite{EK95,BP95}, though not
necessarily with preferred deformations. We have therefore shown that quantum
groups are in fact a special case of star products. For more details see e.g.
the original papers and \cite{FS94,FlS94}.

\subsubsection*{II.4.3 Noncommutative geometry and index theorems}
Noncommutative geometry arose by a kind of ``distillation" from the works
of Connes on $C^*$-algebras and the use in that connection of methods and
results of algebraic geometry. It involves in particular {\it cyclic
cohomology} which was introduced by A. Connes in connection with trace
formulas for operators (cyclic homology was introduced independently by
Tsygan \cite{Ts83}). In particular cyclic cocycles are higher analogues of
traces (see \cite{Di66} for a generalization of the notion of trace). Thus
they facilitate (by setting it algebraically) the computation of the
index, which can obviously be viewed as the trace of some operator, and 
permit to generalize the index theorem, producing {\it algebraic
index formulas} \cite{Co94} of which the Atiyah-Singer formula
(\ref{ASI}) is a special case. As a matter of fact, Fedosov worked first in 
problems related to the index theorem and this brought him naturally to
star product algebras of functions and to the index question in that
context \cite{FeInd} as a fruitful alternative to algebras of
pseudodifferential operators. Recently Nest and Tsygan \cite{NT95} gave 
a nice proof of general algebraic index theorems in the framework
of deformation quantization; doing so they show the existence of a 
``formal trace" (for $X$ symplectic of dimension $2\ell$) given by
${\rm{Tr_{\sc{nt}}}}(u)=\frac{1}{\ell!\nu^{\ell}}\int_X(u\omega^\ell+
\nu\tau_1(u)+\nu^2\tau_2(u)+\cdots)$ where the $\tau_k$ are local expressions
in $u$. That trace satisfies ${\rm{Tr_{\sc{nt}}}}(u*v-v*u)=0$; thus
the integrand will give an equivalence, over ${\Bbb{K}}[\nu^{-1},\nu]]$,
between any given n.c. star product and a strongly closed one.

Cyclic cohomology is based on a bicomplex containing a Hochschild complex
with coboundary operator $b$ of degree 1 and another one with operation $B$
of degree -1 anticommuting with $b$. For a precise definition and properties,
see \cite{Co94}.
The concept does not require to make reference to operator algebras;
formulated abstractly, it applies even better to star products algebras
provided the star products considered are closed (see \cite{CFS92}, where
a explanation of cyclic cohomology in this context can be found).
Indeed, if $*$ is closed (see Def. 2) and a trace $\tau$ is defined on
$u=\sum_{r=0}^{\infty}\nu^r u_r \in N[[\nu]]$ by
$\tau(u)=\int u_{\ell}\omega^{\ell}$, we can consider the quasi-homomorphism
(that measures the noncommutativity of the $*$-algebra and is also a
Hochschild 2-cocycle) $\theta (u_1,u_2) = u_1 * u_2 - u_1u_2$ ; then
${\varphi}_{2k}(u_0,\ldots,u_{2k})=\tau(u_0 *\theta(u_1,u_2) *\ldots*
\theta (u_{2k-1},u_{2k}))$
defines the {\it components of a cyclic cocycle} $\varphi$ in the $(b,B)$
bicomplex on $N$ that is called the {\it character} of the closed star product.
In particular $\varphi_{2\ell}(u_0,\ldots,u_{2 \ell})=
\int u_0 du_1 {\wedge}\ldots {\wedge} du_{2 \ell}$. 
The composition of symbols of pseudodifferential operators is \cite{CFS92}
a closed star product, the character of which coincides with that defined by
the trace on these operators.

A natural extension of the associative algebra context of noncommutative
geometry is to Hopf algebras (in the line of \cite{BFGP94}) and this indeed
permitted now Connes and Moscovici \cite{CM98} to compute the index of
transversally elliptic operators on foliations, a longstanding problem 
(which among many other tools required hypoelliptic pseudodifferential
operators). Another extension, motivated by physics, is to supersymmetric
data, and this has been the subject of recent studies by Fr\"ohlich and
coworkers \cite{FGR98}, first in the context of usual differential geometry
and now in that of noncommutative geometry. There are many more developments
in this framework, including quantized space, but we shall not develop
these further.
\medskip

\noindent {\bf Acknowledgements}. I want to thank Giuseppe (= Joseph) Dito
and Mosh\'e Flato for very useful comments, and Piotr R\c aczka and the
organizers of PFG98 in \L \'od\'z (especially Jakub Rembieli\'nski)
for excellent hospitality in Poland.

\vspace{-2mm}

{\footnotesize
\noindent {\textbf{PACS (1998)}}: 02.10.Tq, 02.10.Vr, 03.65.-w, 02.40.Vh,
11.10.-z

\noindent  {\textbf{MSC (1991)}}: 81S30, 81S10, 81T70, 46M20, 58B30, 58G12,
58F06, 17B37, 19K56, 22E45.}
\end{document}